\documentclass{amsart}
\newcommand{\note}{\noindent {\bf Notation. }}
\newcommand{\remark}{\noindent {\bf Remark. }}
\newcommand{\ws}{\hspace{4pt}}
\newtheorem{theorem}{Theorem}

\newtheorem{corollary}{Corollary}
\newtheorem{lemma}{Lemma}

\usepackage[]{amssymb,amsopn}
\usepackage[]{amsmath}
\usepackage[]{eucal}
\usepackage[]{latexsym}

\begin{document}

\title[Exceptional Laguerre translation]{Translation operator with exceptional Laguerre polynomials}
\author{\'A. P. Horv\'ath }

\thanks{}

\begin{abstract}
We extend the notion of general translation operator to exceptional Laguerre polynomials. To this we investigate the associated singular hyperbolic Cauchy problem.  We derive a maximum principle with respect to this Cauchy problem and applying it we determine the norm of the translation operator. As an application we give Nikol'skii inequalities with respect to exceptional Laguerre polynomials.

\end{abstract}
\maketitle

\section{Introduction}

The classical translation operator $\tau_t$ defined by the formula $(\tau_tf)(x)=f(x+t)$ plays an important role both in approximation theory and harmonic analysis, in particular to characterize and investigate smoothness properties of $f$. The extension of the notion and a general concept is described e.g. in a series of papers by Levitan (\cite{le1}, \cite{le2}, \cite{le3}), in which translation is studied as an integral operator with group-theoretical description, and the underlying space is also examined. The operator is extended to numerical series and to orthogonal polynomials as well.

We concentrate here just one aspect of this rather extended topic. In classical Fourier expansion the kernel function - written in exponential form - can be expressed as $\sum_{k=1}^n\tau_{-t}e^{ikx}$. This ensures that the Dirichlet kernel can be expressed in closed form, and moving the translation from the kernel to the function, convergence theorems can be derived. Christoffel-Darboux formula leads to a similar situation in case of classical orthogonal polynomials, so it is natural to introduce translation on the orthogonal system $\{\varphi_n\}$ as $T_t(\varphi_n(x))=\varphi_n(x)\varphi_n(t)$ and then extend the operator to some function class. Translation operators associated to classical orthogonal polynomials - Jacobi and Laguerre - and to Bessel functions are defined and examined by several authors (see e.g. \cite{ad2}, \cite{m1},  \cite{le4},  \cite{p}). By this operator convolution structures can be defined (cf. e.g. \cite{fjk}, \cite{m1}, \cite{gm}, \cite{cms}), and Nikol'skii-type inequalities can be derived (cf. e.g. \cite{ad}, \cite{ad2}, \cite{adh}). For instance in context of Fourier-Bessel transform the Bessel translation plays exactly the same role as the original translation in Fourier transform  (cf. e.g. \cite{gj}).

This extension can be done by a formal integral operator with the kernel $K(x,s,t)=\sum_{n=0}^{\infty}\varphi_n(x)\varphi_n(s)\varphi_n(t)$ (cf. \cite{m1}). Thus to handle the operator the main tool was to express the kernel in closed form. This led for instance to develop product and addition formulas for orthogonal polynomials (see \cite{a}, \cite{fjk2}, \cite{w}, etc.).

The other method was the investigation of the associated (singular) hyperbolic Cauchy problem cf. e.g. \cite{le4}. By Riemann's method the translation kernel can be obtained in closed form in some cases, for instance the Bessel translation can be expressed as $T^\alpha_t(f;x)={\gamma(\alpha)}\int_0^\pi f\left(\sqrt{t^2+x^2-2xt\cos \varphi}\right)\sin^{2\alpha}\varphi  d\varphi$ (see \cite{le4}) and the Laguerre one as $T_t^\alpha(f;x)=\frac{1}{\sqrt{\pi}}\int_0^\pi f(x+t+2\sqrt{xt}\cos\theta)W^{(\alpha)}(x,t,\theta)d\theta$ (see \cite{w} for the original definition and \cite{adh} for this form). In the applications mentioned above usually the determination of operator norm is the main task. Of course from that nice closed forms it can be derived, but unfortunately Riemann's method can be very difficult if the differential equation is not so simple. The next step was an "error term" method. It was developed by Braaksma and Snoo \cite{b}, \cite{bs} and continued by Markett \cite{m1}. With this method to the final estimation, depending on the coefficients of the differential equation, a simple, double, triple...etc. series has to sum up, so although it works in cases of Bessel and Laguerre translation, but with a bit more complicated differential operator it gives rather inaccurate estimations on the operator norm.

Our aim is to extend the the notion of translation to some kind of exceptional Laguerre polynomials.

The recent development of the area of exceptional polynomial systems is remarkably rapid (see e.g \cite{gukm}, \cite{gumm}, \cite{gfgum}, \cite{lko}, \cite{llms} \cite{ahor} and the references therein). Exceptional polynomials have received contributions both from mathematicians working on orthogonal polynomials and from mathematical physicists. Among the physical applications, exceptional polynomial systems appear mostly as solutions to exactly solvable quantum mechanical problems.

Below we suggest a method similar to the mentioned one, namely to define the translation operator via the associated Cauchy problem, but instead of trying to express the operator in some closed form, we investigate the operator norm directly. This turns our attention to maximum principles with respect to Cauchy problems. This problem has a wide literature for instance in elliptic case, but considering the other cases in hyperbolic case the list of references is rather short. Our main reference points are the  works of Protter and Weinberger \cite{pr}, \cite{prw}, \cite{we}.

This paper is organized as follows. In the second section we derive a maximum principle with respect to a mixed problem on the half-line. In Section 3 we define the translation operator, and by the results of Section 2 we give a method for calculating the operator norm. In Section 4 we introduce exceptional Laguerre polynomials and apply the methods of the previous chapter to define the translation operator and to get the operator norm. Finally as an application of the defined translation operator we prove Nikol'skii type inequalities.

\medskip

\section{Maximum principle}

We examine a special hyperbolic partial differential operator, that is an operator which can be derived from some ordinary (Sturm-Liouville type in the extended sense) differential operator. The symmetry properties of the equation allow to get a maximum principle to the mixed problem generated by the hyperbolic operator.

\subsection{Preliminary results.}
Let
\begin{equation}\label{1} Du=D_xu := u''(x)+q(x)u'(x)-r(x)u(x)=\lambda u(x).\end{equation}
Here $x\in (x_0,x_1)$, and $q\in C^1(x_0,x_1)$, $r\in C(x_0,x_1)$.
Let $u(t)$ be the eigenfunction of $D_tu(t)$ with the same $\lambda$ as above.
As it is usual we take
\begin{equation}\label{ki}\begin{array}{ll}u(t)D_xu(x)=\lambda u(x)u(t)\\
- u(x)D_tu(t)=\lambda u(t)u(x).\end{array}\end{equation}
Denoting  by $u(x,t)=u(x)u(t)$ \eqref{ki} becomes
\begin{equation}\label{75}Lu(x,t)=u_{xx}(x,t)-u_{tt}(x,t)+q(x)u_x(x,t) - q(t)u_t(x,t)-r(x,t)u(x,t)=0,\end{equation}
where
\begin{equation}\label{r}r(x,t)=r(x)-r(t).\end{equation}
So let $L$ be as it is on the left-hand side of \eqref{75}, that is we examine the equation
\begin{equation}\label{35}Lu=u_{xx}-u_{tt}+q(x)u_x - q(t)u_t-r(x,t)u=0,\end{equation}
where subsequently $u=u(x,t)$ is not necessarily in that product form as above.

The formal adjoint operator of $L$ acts on a function $v$ as
$$L^*v=v_{xx}-v_{tt}-q(x)v_x + q(t)v_t-hv,$$
where $$h(x,t)=q'(x)-q'(t)+ r(x,t):=k(x,t)+r(x,t).$$
Thus with $\tilde{u}$ and $v$ (will be given later)
$$vL\tilde{u}-\tilde{u}L^*v=\frac{\partial}{\partial x}\left(v\tilde{u}_x-\tilde{u}v_x+q(x)\tilde{u}v\right)- \frac{\partial}{\partial t}\left(v\tilde{u}_t-\tilde{u}v_t+q(t)\tilde{u}v\right).$$
Our main tool is the Green's formula applying to a characteristic triangle of the differential equation or a subdomain of it.
Since $q$ may have a singularity at $t_0=0$, we define the triangle $\Delta$ as  $\Delta\subset \{(x,t):x \in (0,x_1); 0<t\le x\}$ with vertices $A=(a,m)$; $B=(b,m)$; $C=\left(\frac{a+b}{2},\frac{b-a}{2}+m\right)$ ($0<m\le a< b$). So let $D$ be a subdomain of $\Delta$ such that its border a piecewise smooth curve - possibly with some common part with the edges of $\Delta$. Here we have
\begin{equation}\label{7}\int_{D}vL\tilde{u}-\tilde{u}L^*vd\mu_2=\int_{\partial D}\left(v\tilde{u}_t-\tilde{u}v_t+q(t)\tilde{u}v,v\tilde{u}_x-\tilde{u}v_x+q(x)\tilde{u}v\right)d\mu_1.\end{equation}

Now we begin to define the function $\tilde{u}$.

We take $I=(0,\infty)$. Recalling that in \eqref{35} $q \in C^1(0,\infty)$ and $r \in C[0,\infty)$
with the supplementary condition
\begin{equation}\label{3} u(x,0)= u^0(x) \ws \ws u_t(x,0)=0, \ws \ws x\in I\end{equation}
\eqref{35} defines a Cauchy problem.

In the applications the differential operator $D_x$ has eigenvalues and eigenfunctions $\lambda_k$ and $u_k(x)$ respectively, that is
\begin{equation}\label{4} D_xu_k(x)=\lambda_ku_k(x),\end{equation}
and so with the initial condition
\begin{equation}\label{i} u^0(x)=\sum_{k=0}^n a_ku_k(x)\end{equation}
if $u_k(0)=1$ for all $k\ge 0$, \eqref{35} and \eqref{i} has the solution
\begin{equation}\label{6} u(x,t)=\sum_{k=k_0}^n a_ku_k(x)u_k(t),\end{equation}
of course. First we examine this situation. Furthermore we assume, that $u(x,t)$ is a weighted polynomial, that is
\begin{equation}\label{wp}u(x,t)=p(x,t)w(x,t),\end{equation}
where $p$ is a polynomial, $w$ is a weight function and $w(x,t)>0$ on $(0,\infty)\times (0,\infty)$.

Subsequently we distinguish two cases:

\noindent {\it Case 1.} $r(x,t)\not\equiv 0$. We assume that $u$ is a weighted polynomial.

\noindent {\it Case 2.} $r(x,t)\equiv 0$. $u$ is a solution without any further restrictions.

Let $u$ be a solution of \eqref{35} with initial condition \eqref{3}; in {\it Case 1.} the solution is form of \eqref{6} and the initial condition is like \eqref{i}. Let us choose
\begin{equation}\label{M}\tilde{u}:=|u|-M, \ws  \mbox{in {\it Case 1.}},\ws\ws \tilde{u}:=u-M \ws \mbox{in {\it Case 2.}},\end{equation}
where $M$ is a constant.

Now we define the domains $D$.

\medskip

\begin{lemma} If $u(x,t)$ is a weighted polynomial, then there is a finite partition of $\Delta$ to
$$\Delta=\cup_{i=1}^N\overline{D}_i,$$
where $|u|\in C^2 D_i$ for all $i=1,\dots ,N$; $D_i$ are non-overlapping simply connected domains, and the border of each $D_i$ is continuous, and consists of finitely many smooth arcs: $\gamma_{i,1}\dots \gamma_{i,s_i}$ and all these arcs are on the border of at most two (different) $D_i$-s or it lies on an edge of $\Delta$. That is the system of arcs, domains and points forms a topological image of a polyhedral surface of dimension two. \end{lemma}

\medskip

\proof
Since $w$ is positive on $\Delta$, the sign changes of $u$ are the sign changes of the polynomial, so it is the real algebraic curve $F$. As it is usual, $V(F)$ stands for the zeros, that is the points of the curve itself on the projective plane. Thus, as an application of B\'ezout's theorem (cf. e.g. [Corollary 4.6]\cite{g}) $V(F)$ is  a compact 1-dimensional manifold which means that $V(F)$ is a disjoint union of finitely many connected components, each of which is homeomorphic to a circle (cf. [Remark 5.8]\cite{g}). (Harnack's theorem (cf. e.g. [Proposition 5.10]\cite{g}) gives the upper bound of the loops in case of irreducibile curves, and comparing it to B\'ezout's theorem we have an upper bound of the number of singular points, but we need only the fact that it is finite.) The intersection of this system of arcs and the edges of the triangle results a partition of the triangle which consists of finitely many finitely-connected domains. Inserting finitely many extra arcs carefully (that is inside in a domain in question) one can define the apropriate system of $D_i$-s.

\medskip

We need the following technical lemma.

\begin{lemma}\label{form} Let $u$ be a weighted polynomial solution of \eqref{35} if $r\not\equiv 0$, and an arbitrary solution if $r\equiv 0$. Then with the notation above
$$2(\tilde{u}v)(C)=(\tilde{u}v)(A)+(\tilde{u}v)(B)-\int_{\Delta}-\tilde{u}(v_{xx}-v_{tt}-q(x)v_x+q(t)v_t-kv)+|u|vr d\mu_2 $$ $$+\int_a^b\left(v\tilde{u}_t-\tilde{u}v_t+q(t)\tilde{u}v\right)(x,m)dx$$ $$
+\int_a^{\frac{a+b}{2}}\left(\tilde{u}(2(v_x+v_t)-v(q(t)+q(x)))\right)(x,x-a+m)dx$$ \begin{equation}\label{tri} + \int_{\frac{a+b}{2}}^b\left(\tilde{u}(2(v_t-v_x)-v(q(t)-q(x)))\right)(x,-x+b+m)dx.\end{equation}
\end{lemma}

\medskip

\proof
Let us recall \eqref{7} with $D=D_i$ if $r\not\equiv 0$, and $D=\Delta$ if $r\equiv 0$. Now we take the first case. In the interior of $D_i$ $u$ has constant sign and it is twice continuously differentiable, and by homogeneity $L|u|=0$ here, that is for each $i$
$$\int_{D_i}vL\tilde{u}-\tilde{u}L^*vd\mu_2=\int_{D_i}-\tilde{u}(v_{xx}-v_{tt}-q(x)v_x+q(t)v_t-kv)+|u|vr d\mu_2$$ $$=\int_{\partial D_i}\left(v\tilde{u}_t-\tilde{u}v_t+q(t)\tilde{u}v,v\tilde{u}_x-\tilde{u}v_x+q(x)\tilde{u}v\right)d\mu_1,$$
where the boundary has positive orientation.
Summing up this result with respect to $i$ we have to take into consideration that $\mu_2\left(\cup_i\partial D_i\right)=0$, and because of the orientation of the boundaries of $D_i$-s, the integral on the common boundary of two neighboring domains is taken once with positive and once with negative sign, that is it gives a zero in the sum, that is
$$\int_{\Delta}-\tilde{u}(v_{xx}-v_{tt}-q(x)v_x+q(t)v_t-kv)+|u|vr d\mu_2$$ $$=\sum_{i=1}^N\int_{D_i}-\tilde{u}(v_{xx}-v_{tt}-q(x)v_x+q(t)v_t-kv)+|u|vr d\mu_2$$ $$=\sum_{i=1}^N\sum_{j=1}^{s_i}\int_{\gamma_{i,j}}\left(v\tilde{u}_t-\tilde{u}v_t+q(t)\tilde{u}v,v\tilde{u}_x-\tilde{u}v_x+q(x)\tilde{u}v\right)d\mu_1$$ $$=\sum_{i,j: \gamma_{i,j}\subset \overline{AB}}\int_{\gamma_{i,j}}(\cdot) + \sum_{i,j: \gamma_{i,j}\subset \overline{BC}}\int_{\gamma_{i,j}}(\cdot) +\sum_{i,j: \gamma_{i,j}\subset \overline{CA}}\int_{\gamma_{i,j}}(\cdot)=I+II+III,$$
where $\overline{AB}$ stands for the $AB$ edge of the triangle. The partition of $\Delta$ ensures partitions on the edges.
The first integral of the last line is just
$$I=\int_a^b\left(v\tilde{u}_t-\tilde{u}v_t+q(t)\tilde{u}v\right)(x,m)dx,$$
where $\tilde{u}_t$ is understood piecewise. On the two remaining edges the partition (projected to the $x$ axis) is as it follows $\frac{a+b}{2}=a_1< a_2<\dots < a_K=b$ and $a=b_1< b_2< \dots < b_L=\frac{a+b}{2}$. (As it is usual $a_i$-s are the first coordinates of the partition belonging to side $\overline{BC}$.) On an interval $(a_i,a_{i+1})$ we have
$$II_i:=\int_{a_i}^{a_{i+1}}\left(v(\tilde{u}_x-\tilde{u}_t)-\tilde{u}(v_x-v_t)-\tilde{u}v(q(t)-q(x)))\right)(x,-x+b+m)dx.$$
Let us apply that
$$v(\tilde{u}_x-\tilde{u}_t)-\tilde{u}(v_x-v_t)=\left((\tilde{u}v)(x,-x+b+m)\right)'-2\tilde{u}(v_x-v_t)(x,-x+b+m).$$
(On the right-hand side the notation means that we differentiate the function of one variable $\tilde{u}v(f(x), g(x))$ with respect to its variable.)
Thus we get
$$II_i=(\tilde{u}v)(a_{i+1},-a_{i+1}+b+m)-(\tilde{u}v)(a_{i},-a_{i}+b+m)$$ $$+\int_{a_i}^{a_{i+1}}\left(\tilde{u}(2(v_t-v_x)-v(q(t)-q(x)))\right)(x,-x+b+m)dx.$$
Similarly on $(b_j,b_{j+1})$
$$III_j:=\int_{b_j}^{b_{j+1}}\left(\tilde{u}(v_x+v_t)-v(\tilde{u}_x+\tilde{u}_t)-\tilde{u}v(q(t)+q(x)))\right)(x,-x+b+m)dx,$$
considering that
$$\tilde{u}(v_x+v_t)-v(\tilde{u}_x+\tilde{u}_t)=-\left((\tilde{u}v)(x,x-a+m)\right)'+2\tilde{u}(v_x+v_t)(x,x-a+m)$$
we get
$$III_j=-(\tilde{u}v)(b_{j+1},b_{j+1}-a+m)+(\tilde{u}v)(b_{j},b_{j}-a+m)$$ $$+\int_{b_j}^{b_{j+1}}\left(\tilde{u}(2(v_t+v_x)-v(q(t)+q(x)))\right)(x,-x+b+m)dx.$$
Thus
$$\sum_{i=1}^{K-1}II_i+\sum_{j=1}^{L-1}III_j=(\tilde{u}v)(A)+(\tilde{u}v)(B)-2(\tilde{u}v)(C)$$ $$+\int_a^{\frac{a+b}{2}}\left(\tilde{u}(2(v_x+v_t)-v(q(t)+q(x)))\right)(x,x-a+m)dx $$ $$+ \int_{\frac{a+b}{2}}^b\left(\tilde{u}(2(v_t-v_x)-v(q(t)-q(x)))\right)(x,-x+b+m)dx.$$
After a rearrangement these computations give the statement of the lemma in the first case.

When $r\equiv 0$, $D=\Delta$ and $\tilde{u}=u-M$. Following the computation above we arrive the statement again.

\medskip

\subsection{Cauchy problem, maximum principle.}

Below we give a maximum principle for hyperbolic Cauchy problem. When $r(x,t)\not\equiv 0$ the statement is restricted to weighted polynomial solutions, if there are any. For simplicity let us denote by $$L^1v:=v_{xx}-v_{tt}-q(x)v_x+q(t)v_t-kv.$$

\begin{theorem}\label{T1} Let  $r(x,t)\ge 0$ on $x>t>0$ and let us assume that $u(x,t)=\sum_{k=0}^n a_ku_k(x)u_k(t)$ is the weighted polynomial solution to the Cauchy problem \eqref{35} with initial condition \eqref{3}, where $u^0(x)=\sum_{k=0}^n a_ku_k(x)$; or $r(x,t)\equiv 0$  on $x>t>0$ and $u(x,t)$ is a solution to \eqref{35} and \eqref{3}. Additionally  we assume that and $\lim_{\|(x,t)\|\to \infty}u(x,t)=0$. Let us assume that there is a function $v(x,t)$ such that
$$ v>0 \ws \ws L^1v> 0 \ws \ws \mbox{on} \ws x>t>0;$$ $$  2(v_t+v_x)-v(q(t)+q(x))>0 \ws \ws \mbox{on} \ws t=x-K,\ws \forall K>0;$$ $$  \ws \ws 2(v_t-v_x)-v(q(t)-q(x))>0 \ws \ws \mbox{on} \ws t=-x+L,\ws \forall L>0$$
$$q(t)v-v_t>0 \ws \ws \mbox{on} \ws t=m, \ws \forall m<m_0$$
with some $m_0$ small enough.
Then
$$\|u(x,t)\|_{\infty, (0\infty)\times (0,\infty)}=\|u^0(x)\|_{\infty,(0,\infty)}.$$\end{theorem}

\medskip

\proof
Let $\|u^0(x)\|=M_1$ and $\varepsilon>0$ arbitrary. Let $M:=M_1+\varepsilon$. We define $\tilde{u}$ as in \eqref{M}. Thus $\tilde{u}\le-\varepsilon$ on the halfline and at infinity. Because of symmetry it is enough to deal with the lower triangle $t\le x$. If $\tilde{u}$ is positive somewhere on the lower triangle, then there is a point $C=(c_1,c_2)$, such that $\tilde{u}(C)=0$ while $\tilde{u}(x,t)<0$ on the characteristic triangle with vertex $C$. (We can choose $C$ as $c_2=min\{\tilde{c_2}:\tilde{u}(\tilde{c_1},\tilde{c_2})\ge 0\}$ This $c_2$ is not unique, but $0<c_2\le c_1<\infty$, because of the assumption on the initial line and at infinity.)  Let $c_2>m>0$ arbitrary small. We define $\Delta$ as above with this $C$. Recalling Lemma 2, by the assumptions on $v$ it can be seen that the left-hand side, $2(\tilde{u}v)(C)=0$ and oppositely on the right-hand side $(\tilde{u}v)(A)$ and $(\tilde{u}v)(B)$ are negative by the definition of $C$, the integrand on the triangle is positive and on the edges $\overline{BC}$ and $\overline{CA}$  are negative. On $\overline{AB}$ the initial condition ensures that the first term tends to zero, when $m$ tends to zero, and the other terms are also negative. That is, if $m$ is small enough, on the left-hand side each term is negative which is a contradiction.

It directly gives the statement when $u$ is a weighted polynomial. If $r(x,t)\equiv 0$, we have to repeat the previous chain of ideas with $-u$ and these results together give the statement again.

\medskip

\remark

(1) If $r(x,t)\equiv 0$ the proof shows that the operator which maps $u^0$ to $u$ is positive.

(2) As a corollary of Theorem 1 one can derive more general maximum principles by extension from dense subspace. We discuss it by an example later.

\medskip

\section{Translation operator}

\subsection{Translation operator in finite dimension}
Let $w$ be a weight function on $[x_0,x_1)$ and $\{\varphi_k\}_{k=0}^{\infty}$ an orthogonal system in $L^2_w$, that is
$$\int_{x_0}^{x_1}\varphi_k(x)\varphi_n(x)w(x)dx=\delta_{k,n}\sigma^2_n, \ws\ws k,n=0,1,\dots .$$
Thus denoting by
\begin{equation}\label{22}\mathcal{P}_n:=\mathrm{span}\{\varphi_k\}_{k=0}^n; \ws \ws \mathcal{P}:=\cup_n\mathcal{P}_n,\end{equation}
for all $t\in [x_0,x_1)$ an operator $T_t$ acting on $\mathcal{P}$ can be defined as
\begin{equation}\label{23}T_t\left(\sum_{k=0}^n a_k\varphi_k(x)\right)=:T_t(u,x)=\sum_{k=0}^n a_k\varphi_k(x)\varphi_k(t).\end{equation}
The operator norm in the space $L^2_w$ can be investigated. We introduce the following notation.

\medskip

\note $\|T_t\|_{p,w}$ means the operator norm when the operator acts on  $L^p_w$ and maps into  $L^p_w$.

\medskip

Let us suppose that
\begin{equation}\label{n1}\|\varphi_k\|_{\infty,[x_0,x_1)}=1 \ws \ws \forall \ws k\in\mathbb{N}.\end{equation}
Since
$$\left\|\sum_{k=0}^n a_k\varphi_k(x)\right\|_{2,w}^2=\sum_{k=0}^n (a_k\sigma_k)^2,$$
and
$$\left\|\sum_{k=0}^n a_k\varphi_k(t)\varphi_k(x)\right\|_{2,w}^2=\sum_{k=0}^n (a_k\varphi_k(t)\sigma_k)^2\le \sum_{k=0}^n (a_k\sigma_k)^2,$$
the operator norm in $L^2_w$ is obviously at most one. If there is a $t$ and a $k$ with $\varphi_k(t)=1$, considering $T_{t}\varphi_k(x)=\varphi_k(x)$, we have that $\|T_{t}\|_{2,w}=1$.

\medskip

\remark
We use the notation $[x_0,x_1)$, because in the next section we concentrate to the half-line. The same method can be applied to a finite interval, see for instance the Jacobi translation (cf. \cite{ad2}).

\medskip

\subsection{Translation operator by Cauchy problem}

Returning to our Cauchy problem, if it is known that our singular mixed problem \eqref{35}, \eqref{3} with an additional assumption at infinity has a unique symmetric solution, we can define a linear operator $T$ on $L_{[0,\infty),\infty}$ as it follows.
\begin{equation}\label{46} T (u^0(x))=u(x,t)=(Tu^0)(x,t).\end{equation}

To make some comparison between the previous operator and the finite dimensional case, we are interested in (symmetric) solutions on $[0,\infty)$ as \eqref{6} of \eqref{35}, with initial condition \eqref{3} in form \eqref{i}.

Let us recall that $u_k$-s are the eigenfunctions of the differential operator \eqref{4}. If $\{u_k\}_{k=0}^{\infty}$ fulfils the conditions
\begin{equation}\label{36} \lim_{x\to 0+ \atop x\to \infty} u_k(x)u'_m(x)e^{\frac{1}{2}\int^xq}=0,\ws\ws k,m \in \mathbb{N}\end{equation}
\begin{equation}\label{47} u_k(0)=1, \ws \ws \ws u'_k(0)=0 \ws \ws \forall \ws k=0,1,\dots, \end{equation}
\begin{equation}\label{48} \|u_k\|_{\infty}=1 \ws \ws \forall  \ws k\in\mathbb{N},\end{equation}
then by \eqref{36} the eigenfunctions form an orthogonal system on $(0,\infty)$ with respect to the weight $w(x)=e^{\int^xq}$ and by \eqref{47}
the solutions of the form  \eqref{6} obviously fulfils the initial condition \eqref{i}. In addition by \eqref{48} we can apply the chain of ideas with respect to $\{\varphi_k\}_{k=0}^{\infty}$ to $\{u_k\}_{k=0}^{\infty}$. Thus if the eigenfunctions of \eqref{4} satisfy the conditions \eqref{36} and \eqref{47}, we will obtain that if $u(x)\in \mathcal{P}$,
\begin{equation}\label{49}(T u)(x,t)=T_t (u,x).\end{equation}
That is $T: \mathcal{P}\to C_{[0,\infty)\times [0,\infty)}$ and for all $t\in I$ $T_t: \mathcal{P} \to C_{[0,\infty)}$ and
\begin{equation}\label{60} \|T|_\mathcal{P}\|_{\infty}=\sup_{t\in I}\|T_t|_\mathcal{P}\|_{\infty}.\end{equation}
That is if there is a maximum principle similar to Theorem \ref{T1} it will show that similarly to the $L^2_w$ case, for all $t\in [0,\infty)$, $\|T_t\|_{\infty}\le 1$ on $\mathcal{P}$.

If we did not have an existence and uniqueness theorem with respect to the singular mixed problem above or a general maximum principle is not known, according to the previous chain of ideas we can follow the next process.
If $\mathcal{P}$ is dense in $L_{[0,\infty),\infty}$, we can extend the operator from $\mathcal{P}$ to the whole space with the same operator norm and so we can define the translation operator with formula \eqref{49} again.

Let us observe that in general the orthogonality of the eigenfunctions ensures some "selfadjointness" of the operator on $\mathcal{P}$:
$$\int_{0}^{\infty}T_t(u_k,x)u_l(x)w(x)dx=\int_{0}^{\infty}u_k(x)u_k(t)u_l(x)w(x)dx=u_k(t)\delta_{k,l}\sigma_k^2$$ \begin{equation}\label{s}=\int_{0}^{\infty}u_k(x)T_t(u_l,x)w(x)dx.\end{equation}
Considering the linearity the same fulfils for any $p(x)\in \mathcal{P}_n$ and $q(x)\in \mathcal{P}_m$. Now let us assume, that $\mathcal{P}$ is dense in $L_{\infty}$ (on $[0,\infty)$). Since (with a weight function $w$) for any $p(x)\in \mathcal{P}$
$$\|T_t(p,x)\|_{1,w}=\sup_{q\in \mathcal{P} \atop \|q\|_{\infty}=1}\int_{0}^{\infty}T_t(p,x)q(x)w(x)dx=\sup_{q\in \mathcal{P} \atop \|q\|_{\infty}=1}\int_{0}^{\infty}p(x)T_t(q,x)w(x)dx$$ \begin{equation}\label{d}\le \sup_{q\in \mathcal{P} \atop \|q\|_{\infty}=1}\|T_t(q,x)\|_{\infty}\|p\|_{1,w}\le \|p\|_{1,w},\end{equation}
we have that for all  $t\in [0,\infty)$ $\|T_t|_{\mathcal{P}}\|_{1,w}\le 1$. If $\mathcal{P}$ is dense in $L_{p,w}$, $1\le p \le \infty$ to estimate the operator norm in between $L_{p,w}$- spaces we can apply the Riesz-Thorin theorem. We give a detailed description of the method in the next section.

\medskip

\subsection{Examples} Below we give examples for application of Theorem \ref{T1}. Here we deal with the special potential
$$q(x)=\frac{2\alpha+1}{x} >0, \ws \ws x\in(0,x_1).$$
It means that
\begin{equation}\label{24}Lu=u_{xx}-u_{tt}+\frac{2\alpha+1}{x}u_x - \frac{2\alpha+1}{t}u_t-ru.\end{equation}
By different choices of $r$ we get the Bessel, Laguerre and exceptional Laguerre equations and translations. In the next section we deal with exceptional Laguerre translation.

\medskip

 The normalized Bessel functions are
\begin{equation}\label{26}j_{\alpha}(z)=\Gamma(\alpha+1) \left(\frac{2}{z}\right)^\alpha J_\alpha(z)= \sum_{k=0}^\infty\frac{(-1)^k\Gamma(\alpha+1)}{\Gamma(k+1)\Gamma(k+\alpha+1)}\left(\frac{z}{2}\right)^{2k}.\end{equation}
 $j_\alpha(xy)$ are the solution of the problem  $u''(x)+\frac{2\alpha+1}{x}u'(x)+y^2u(x)=0$,  $u(0)=1$, $u'(0)=0$, that is we can derive the Cauchy problem \eqref{35} and \eqref{3} with
\begin{equation}\label{25}r(x,t)\equiv 0.\end{equation}
One can take into consideration Bessel translation on the half-line, see e.g. \cite{le4}.

\medskip

The Laguerre differential equation is
\begin{equation}\label{27}xy''(x)+\left(\alpha+1-x\right)y'(x)+ny(x)=0.\end{equation}
That is
\begin{equation}\label{28}u_n(x):=\frac{n!\Gamma(\alpha+1)}{\Gamma(n+\alpha+1)}L_n^{(\alpha)}(x^2)e^{-\frac{x^2}{2}}\end{equation}
fulfils the eigenfunction equation
\begin{equation}\label{29}u''(x)+\frac{2\alpha+1}{x}u'(x)-x^2u(x)=-4\left(n+\frac{\alpha+1}{2}\right)u(x).\end{equation}
Thus
\begin{equation}\label{30}r(x,t)=x^2-t^2.\end{equation}
As above, since $u_n(0)=1$ for all $n\in\mathbb{N}$ thus a solution like \eqref{6} fulfils the initial condition \eqref{i}. The Laguerre translation is discussed e.g. in \cite{gm} and \cite{adh}.

\medskip

Now we derive the known result for the Bessel and Laguerre translation from Theorem \ref{T1}. Let $u_k(x)=j_\alpha(\lambda_k x)$, where $\lambda_k$-s are the positive zeros of Bessel functions; or $u_k(x)$ is given by \eqref{28}. Then

\begin{theorem}\label{TN} Let us consider the Cauchy problem \eqref{35} with $q(x)=\frac{2\alpha+1}{x}$ and $r(x,t)=0$ with initial condition \eqref{3}; or $r(x,t)=x^2-t^2$; with initial condition \eqref{3}, where $u^0(x)=\sum_{k=0}^n a_ku_k(x)$. Then the solution $u(x,t)$; or $u(x,t)=\sum_{k=0}^n a_ku_k(x)u_k(t)$ respectively fulfils the maximum principle
$$\|u(x,t)\|_{\infty, (0\infty)\times (0,\infty)}=\|u^0(x)\|_{\infty,(0,\infty)}.$$\end{theorem}

\medskip

\proof According to Theorem \ref{T1} it is enough to define an appropriate function $v$. So let
\begin{equation}\label{v} v(x,t):=x^{1+\alpha}t^{1+\alpha}.\end{equation}
By this $v=v(x,t)$ we can compute
$$2(v_t\pm v_x)-v(q(t)\pm q(x))$$\begin{equation}\label{vs}=v\left(2(1+\alpha)\left(\frac{1}{t}\pm \frac{1}{x}\right)-(2\alpha+1)\left(\frac{1}{t}\pm \frac{1}{x}\right)\right)=v\left(\frac{1}{t}\pm \frac{1}{x}\right).\end{equation}
Similarly we have
$$v_{xx}-v_{tt}-q(x)v_x+q(t)v_t-kv$$ $$=v\left(-\alpha(1+\alpha)\left(\frac{1}{t^2}- \frac{1}{x^2}\right)+(2\alpha+1)(1+\alpha)\left(\frac{1}{t^2}- \frac{1}{x^2}\right)-(2\alpha+1)\left(\frac{1}{t^2}- \frac{1}{x^2}\right)\right)$$
\begin{equation}\label{vt}=\alpha^2v\left(\frac{1}{t^2}- \frac{1}{x^2}\right).\end{equation}
Finally on $x\in(a,b); t=m$ we need
\begin{equation}\label{f} v\left(q(t)-\frac{v_t}{v}\right)=v\frac{\alpha}{t}.\end{equation}
Recalling that we are on the lower triangle, the right-hand sides of \eqref{vs}, \eqref{vt}  and \eqref{f}  are positive which proves the statement.

\medskip

\remark

\noindent (1) In  Laguerre case we follow the the proof of Theorem 1 via Lemma 1 and Lemma 2. It gives that the norm of the (restricted) operator is bounded by 1, and the proof does not ensure the positivity of the operator. Indeed Laguerre translation is not a positive operator, cf. \cite{gm}. In Bessel case, taking into consideration the first remark after Theorem \ref{T1}, we get that Bessel translation is a positive operator, cf. \cite{le4}.

\noindent (2) For $J_n(\lambda_kz)$ Bessel case see \cite{pr}.

\medskip

\section{Exceptional Laguerre translation and Nikol'skii inequality}

\subsection{Exceptional Laguerre polynomials}  Exceptional polynomials of codimension $m\geq 1$ are real-valued polynomial sequences, where $m$ degrees are missing from the sequence of degrees. All the polynomials in the sequence are eigenpolynomials of a second order real valued linear differential operator with rational coefficients. Furthermore exceptional polynomials form an orthogonal system with respect to a weight function on a real interval, and in the weighted Hilbert space in question this system is closed (cf. \cite[Definition 7.4]{gfgum}). As it is pointed out in the previously cited work, exceptional operators can be obtained by finite many Crum-Darboux transformations from Bochner operators.

We are interested in exceptional Laguerre polynomials on the half-line. Applying one-step Darboux transformation three different types of exceptional Laguerre polynomials appear in the literature (cf. \cite{gumm} and \cite{llms}), $\{L_{m,m+n}^{i, (\alpha)}\}_{n=0}^{\infty}$ for $i=I$ and $i=II$ and $\left(1\cup\{L_{m,m+n}^{i, (\alpha)}\}_{n=1}^{\infty}\right)$ for $i=III$. (The second index shows the degree of the polynomial.)

In case $I$ $\alpha>0$, in case $II$ $\alpha>m-1$ and in case $III$ $-1<\alpha<0$.\\ $\{L_{m,m+n}^{i, (\alpha)}\}$ are the orthogonal polynomials on $(0, \infty)$ with respect to the weight
$$w:=w^{(\alpha)}_m:=\frac{x^{\alpha}e^{-x}}{S_i^2(x)},$$
where $S_i=S^{(\alpha)}_{i,m}$
$$S_I(x):=L_m^{(\alpha-1)}(-x), \ws \ws \ws S_{II}(x):=L_m^{(-\alpha-1)}(x) \ws \ws \ws S_{III}(x):=L_m^{(-\alpha-1)}(-x).$$  In this subsection we investigate these three types simultaneously.

$L_{m,m+n}^{i, (\alpha)}$ satisfies the differential equation (cf. \cite[(8), (13)]{ahor}, \cite[5.13]{llms})
\begin{equation}\label{31}xy''(x)+\left(\alpha+1-x-2x\frac{S'_i(x)}{S_i(x)}\right)y'(x)+\left(n+\varepsilon_im-2\delta_i\alpha\frac{S'_i(x)}{S_i(x)}\right)y(x)=0,\end{equation}
where $\varepsilon_I=1$, $\varepsilon_{II}=-1$, $\varepsilon_{III}=0$ $\delta_I=\delta_{II}=1$, $\delta_{III}=0$.
That is
\begin{equation}\label{32}u_{n,i}(x):=c_{n,i}L_{m,m+n}^{i, (\alpha)}(x^2)\frac{e^{-\frac{x^2}{2}}}{S_i(x^2)}\end{equation}
fulfils the eigenfunction equation
$$u''(x)+\frac{2\alpha+1}{x}u'(x)$$ $$-\left(x^2 +4(\alpha+1-2\delta_i+x^2)\frac{S_i'}{S_i}(x^2)-4x^2\frac{S''_i}{S_i}(x^2)+8x^2\left(\frac{S'_i}{S_i}\right)^2(x^2)\right)u(x)$$ \begin{equation}\label{33}=-4\left(n+\varepsilon_im+\frac{\alpha+1}{2}\right)u(x).\end{equation}
That is
$$r_i(x,t)=x^2-t^2+\left(4(\alpha+1-2\delta_i+x^2)\frac{S'_i}{S_i}(x^2)-4x^2\frac{S''_i}{S_i}(x^2)+8x^2\left(\frac{S'_i}{S_i}\right)^2(x^2)\right.$$ \begin{equation}\label{34}\left.-4(\alpha+1-2\delta_i+t^2)\frac{S'_i}{S_i}(t^2)-4t^2\frac{S''_i}{S_i}(t^2)+8t^2\left(\frac{S'_i}{S_i}\right)^2(t^2)\right).\end{equation}
Furthermore $\{u_{n,i}(x)\}$ is a sequence of orthogonal functions with respect to the weight $x^{2\alpha+1}$ that is $\int_0^{\infty}u_{n,i}(x)u_{m,i}(x)x^{2\alpha+1}dx=\sigma_{n,i}^2\delta_{m,n}$. An appropriate choice of $c_{n,i}$ ensures that  $u_{n,i}(0)=1$, and the solution of \eqref{35} in form like  \eqref{6} fulfils the initial condition \eqref{i}.

\subsection{The uniform norm on the half-line} In this subsection we deal with exceptional Laguerre polynomials of type I and II. The next lemma is about $I$-type exceptional Laguerre polynomials. It ensures that a weighted polynomial $u(x,t)=\sum_{k=1}^na_ku_k(x)u_k(t)$ fulfils the initial condition $u(x,0)=\sum_{k=1}^na_ku_k(x)$, and simultaneously the operator norm (at least in $L^2$) is one.

\medskip

\begin{lemma}\label{Ll1} Let $$z_n(x)=L_{m,m+n}^{I, (\alpha)}\frac{e^{-\frac{x}{2}}}{L_m^{(\alpha-1)}(-x)}.$$ Then
\begin{equation}\label{50}\|c_nz_{n}\|_{\infty}=u_n(0)=1\ws\ws n=0,1 \dots.\end{equation}\end{lemma}

\medskip

\proof
The derivation of exceptional polynomials from the classical ones ensure the following formula \cite[(21)]{gumm}.
$$\frac{L_{m,m+n}^{I, (\alpha)}(x)}{L_m^{(\alpha-1)}(-x)}=\frac{L_m^{(\alpha)}(-x)L_{n}^{ (\alpha)}(x)-L_{m-1}^{(\alpha)}(-x)L_{n-1}^{ (\alpha)}(x)}{L_m^{(\alpha-1)}(-x)}.$$
By this one and by \cite[(5.1.13)]{sz} we have
$$\frac{L_{m,m+n}^{I, (\alpha)}(x)}{S_I(x)}=\frac{L_{n}^{ (\alpha)}(x)L_m^{(\alpha-1)}(-x)+L_{m-1}^{(\alpha)}(-x)L_n^{(\alpha-1)}(x)}{L_m^{(\alpha-1)}(-x)}$$ \begin{equation}\label{51}=L_n^{(\alpha)}(x)+\frac{L_{m-1}^{(\alpha)}(-x)}{L_m^{(\alpha-1)}(-x)}L_n^{(\alpha-1)}(x).\end{equation}
Since $e^{-\frac{x}{2}}L_n^{(\alpha)}(x)$ attains its $\sup$-norm at zero for all $n$ (see e.g. \cite[(2.6)]{adh}), it is enough to show that the same is valid for
$\frac{L_{m-1}^{(\alpha)}(-x)}{L_m^{(\alpha-1)}(-x)}$. By \cite[(5.1.14)]{sz} and \cite[(5.1.6)]{sz} we have
$$ \left(L_m^{(\alpha-1)}(-x)\right)^2\left(\frac{L_{m-1}^{(\alpha)}(-x)}{L_m^{(\alpha-1)}(-x)}\right)'=\sum_{k=1}^{m-1}\binom{m-1+\alpha}{m-1-k}\frac{x^{k-1}}{(k-1)!}\sum_{j=0}^m\binom{m+\alpha-1}{m-j}\frac{x^{j}}{j!}$$ $$-\sum_{k=0}^{m-1}\binom{m-1+\alpha}{m-1-k}\frac{x^{k}}{k!}\sum_{j=1}^m\binom{m+\alpha-1}{m-j}\frac{x^{j-1}}{(j-1)!}$$ $$=\sum_{1\le k, j\le m-1}\binom{m-1+\alpha}{m-1-k}\binom{m+\alpha-1}{m-j}\frac{1}{(k-1)!(j-1)!}\left(\frac{1}{j}-\frac{1}{k}\right)x^{k+j-1}$$ $$+\binom{m+\alpha-1}{m}\sum_{k=1}^{m-1}\binom{m-1+\alpha}{m-1-k}\frac{x^{k-1}}{(k-1)!}$$ $$-\binom{m+\alpha-1}{m-1}\sum_{j=1}^{m}\binom{m-1+\alpha}{m-j}\frac{x^{j-1}}{(j-1)!}$$ $$+\frac{x^m}{m!}\sum_{k=1}^{m-1}\binom{m-1+\alpha}{m-1-k}\frac{x^{k-1}}{(k-1)!}-\frac{x^{m-1}}{(m-1)!}\sum_{k=0}^{m-1}\binom{m-1+\alpha}{m-1-k}\frac{x^{k}}{k!}$$ $$=\Sigma_1+\Sigma_2+\Sigma_3.$$
$$\Sigma_2=-\frac{\alpha+m}{m}\sum_{k=1}^{m-1}\binom{m+\alpha-1}{m-1}\binom{m-1+\alpha}{m-1-k}\frac{k}{m-k}\frac{x^{k-1}}{(k-1)!}$$ $$-\binom{m+\alpha-1}{m-1}\frac{x^{m-1}}{(m-1)!}<0.$$
$$\Sigma_3=\sum_{k=1}^{m-1}\binom{m-1+\alpha}{m-1-k}\frac{1}{(k-1)!(m-1)!}\left(\frac{1}{m}-\frac{1}{k}\right)x^{m+k-1}$$ $$-\binom{m-1+\alpha}{m-1}\frac{x^{m-1}}{(m-1)!}<0.$$
Since $\Sigma_1=\sum_{1\le k, j\le m-1}a(k,j)x^{k+j-1}$, where $a(k,k)=0$ and for $j\neq k$ we compute
$$a(k,j)+a(j,k)$$ $$=\binom{m-1+\alpha}{m-k}\binom{m-1+\alpha}{m-j}\frac{1}{(k-1)!(j-1)!}(m+\alpha)\left(\frac{1}{j}-\frac{1}{k}\right)(j-k)<0,$$
that is
$$\Sigma_1<0.$$
Thus $\frac{L_{m-1}^{(\alpha)}(-x)}{L_m^{(\alpha-1)}(-x)}$ is decreasing on $[0,\infty)$ and consequently $R_{m,m+n}^{I, (\alpha)}(x)\frac{e^{-\frac{x}{2}}}{L_m^{(\alpha-1)}(-x)}$ attains its maximum at zero, where
\begin{equation}\label{n} R_{m,m+n}^{I, (\alpha)}(x)=\frac{L_{m,m+n}^{I, (\alpha)}(x)}{\binom{n+\alpha-1}{n-1}\frac{n+\alpha+m}{n}}.\end{equation}
Since $\left|L_{n}^{ (\alpha)}(x)\right|e^{-\frac{x}{2}}$ attains its maximum only at zero, this fulfils also for $z_n(x)$.

\medskip

Previous computation used the formula which expresses exceptional polynomials by the classical ones and the corresponding property of the classical polynomials. Actually all the information have to be contained by the differential equation.
Below we adapt a method to derive an estimation to the norm of $z_n(x)$ by $z_n(0)$ by the differential equation. In the next lemma we apply this idea to $II$-type exceptional Laguerre polynomials.

Our main tool is the following ( cf.\cite[§ 7.31, Theorem 7.6.1]{sz}). On an interval $(a,b)$ let $\varphi$ be a solution of a differential equation:
\begin{equation}\label{65}(kz')'+\Phi_1 z=0.\end{equation}
Let $f=z^2+\frac{k}{\Phi_1}(z')^2$. Then
$$f'=-(k\Phi_1)'\frac{(z')^2}{\Phi_1^2}.$$
Let us observe that $f=z^2$ if $z'=0$ (or $\frac{k}{\Phi_1}=0$).
If $k>0$, $\Phi_1 >0$ for $a<x<z_0$ (or possibly $\frac{k}{\Phi_1}(a)=0$) then the sequence formed by ($|\varphi(a)|$) and by the relative maxima of $|\varphi|$ is strictly decreasing when $(k\Phi_1)'>0$ and is strictly increasing when $(k\Phi_1)'<0$.

To prove a result like \eqref{50} it is enough to show that $z_{n,i}(x):=L_{m,m+n}^{i, (\alpha)}(x)\frac{e^{-\frac{x}{2}}}{S_i(x)}$ attains its maximum at zero. $z_{n,i}$ is the solution of the differential equation
$$\left(x^{\alpha+1}z'\right)'$$ $$+x^{\alpha+1}\left(\frac{n+_i m+\frac{\alpha+1}{2}}{x}+\left(\frac{1+(1-2\delta_i)\alpha}{x}-1\right)\frac{S_i'}{S_i}-2\left(\frac{S'_i}{S_i}\right)^2+\frac{S''_i}{S_i}-\frac{1}{4}\right)z$$ \begin{equation}= \left(x^{\alpha+1}z'\right)'+x^{\alpha+1}\Phi z=0,\end{equation}
where $\Phi=\Phi_{n,m,i,\alpha}$.

That is it is enough to show that there is an $A$ such that $2(\alpha+1)\Phi(x)+x\Phi'(x)>0$ on $(0,A)$, and $|z_n(x)|<z_n(0)$ if $x\ge A$.

\medskip

\begin{lemma}\label{Ll2} Let
$$z_n(x)=L_{1,1+n}^{II, (\alpha)}(x)\frac{e^{-\frac{x}{2}}}{S(x)},$$
where  $\alpha\ge 1$, $S(x)=-x-\alpha$. Then
\begin{equation}\|c_nz_{n}\|_{\infty}=u_n(0)=1\ws\ws n=0,1 \dots.\end{equation}\end{lemma}

\medskip

\proof
According to \cite[(35)]{gumm}
$$z_n(x)=e^{-\frac{x}{2}}\left(-xL_{n-1}^{\alpha+2}(x)+\alpha\left(1+\frac{1}{x+\alpha}\right)L_{n}^{\alpha+1}(x)\right).$$
Considering \cite[Theorem 8.22.8]{sz} $A=2n+\alpha+1$ is an appropriate choice.
$$2(\alpha+1)\Phi(x)+x\Phi'(x)$$ $$=(2\alpha+1)n+\alpha^2-\frac{5}{2}\alpha-\frac{1}{2}-\frac{\alpha+1}{2}x-2\frac{\alpha+2}{\alpha+x}+\frac{\alpha(4\alpha+5}{(\alpha+x)^2}+\frac{4x^2}{(\alpha+x)^3}>0,$$
when $n\ge 3$ an $0<x<2n+\alpha+1$. For $n=0,1,2$ one can easily check that the norm is attained at zero.

\medskip

\remark
By this method can be proved that the weighted Jacobi and Laguerre polynomials attain their maximum at the endpoint of the interval in question.

\subsection{Density} Subsequently we study the exceptional Laguerre polynomials of type I. Let $S=S_I$. In order to apply the second definition of translation and to estimate the operator norm in $L^q_w$ the density question has to be investigated. Our aim is to show that $\mathrm{span}\{L_{m,m+n}^{I, (\alpha)}(x)\frac{e^{-\frac{x}{2}}}{S(x)}\}$ is dense in the space $C_v$ and in $L^q_w$ (see below).
$$C_v:=\left\{f\in C(0,\infty): \lim_{x\to \infty} f(x)v=0\right\},$$
where
$$v(x):=w_I(x)=\frac{e^{-\frac{x}{2}}}{S(x)}.$$
Let $\alpha>0$, $w(x)=x^{\alpha}$. Let us recall the notation that for $1\le q <\infty$ the norm in $L^q_{w}$ and the inner product are
$$\|f\|_{q,w}^q=\int_0^{\infty}|f|^qw, \ws \ws \langle f,g\rangle =\int_0^{\infty}fgw$$
respectively.
Let us define the real linear spaces
$$\mathcal{L}^{I, (\alpha)}_n:=\mathrm{span}\left\{L_{m,m+k}^{I, (\alpha)}(x)\right\}_{k=0}^n, \ws \ws \mathcal{L}^{I, (\alpha)}=\cup_n\mathcal{L}^{I, (\alpha)}_n.$$
First we observe that $S\in \Pi_m$ and it has $m$ different zeros $x_m< \dots <x_1<0$, cf. \cite[5.1.6]{sz}, where $\Pi_k$ is the set of polynomials of degree at most $k$.
According to \eqref{31} at the zeros of $S$ for all $n\ge 0$ $L_{m,m+n}^{I, (\alpha)}$ satisfies the equation
$$x_jp'(x_j)+\alpha p(x_j)=0.$$
Let us define
$$\mathcal{M}^{(\alpha)}_n:=\left\{p\in \Pi_{n+m} : x_jp'(x_j)+\alpha p(x_j)=0, \ws \ws i=1, \dots, m.\right\}.$$
Let $\Pi=\cup_k\Pi_k$ and $\mathcal{M}^{(\alpha)}:=\cup_n\mathcal{M}^{(\alpha)}_n$. In \cite[Lemma 3.1]{lko} it is pointed out that
$$\mathcal{L}^{I, (\alpha)}_n=\mathcal{M}^{(\alpha)}_n.$$
Similarly to \cite[Lemmas 3.1 and 3.2]{gukm} we show that

\begin{lemma}\label{L3} $\mathcal{L}^{I, (\alpha)}v$ is dense in $C_v$ and in $L^q_w$, where $1\le q < \infty$ .\end{lemma}

\proof
Let $f\in C_v$. Then $\frac{f}{S^2}\in C_{\hat{v}}$, where $\hat{v}(x)=S(x)e^{-\frac{x}{2}}$. Obviously the polynomials are dense in $C_{\hat{v}}$. Thus for all $\varepsilon >0$  there is a  polynomial $p$ such that
$$\varepsilon >\left\|\left(\frac{f}{S^2}-p\right)\hat{v}\right\|_{\infty}= \left\|\left(f-pS^2\right)v\right\|_{\infty}.$$
Since $pS^2 \in \mathcal{M}^{(\alpha)}$, the first part of the lemma is proved.

The second part can be proved similarly. $f \in L^q_w$ if and only if $\frac{f}{S} \in L^q_{S^qw}$. The set $Q:=\{p(x)e^{-\frac{x}{2}}: p\in \Pi\}$ is dense in $L^q_{S^qw}$.  Thus for all $\varepsilon >0$ there is a $p \in Q$ such that
$$\varepsilon >\int_0^{\infty}\left|\frac{f}{S(x)}-p(x)e^{-\frac{x}{2}}\right|^qS^q(x)x^{\alpha}dx=\int_0^{\infty}\left|f(x)-(pS^2)(x)v(x)\right|^qS^q(x)x^{\alpha}dx,$$
and we finish the proof again with $pS^2 \in \mathcal{M}^{(\alpha)}$.

\medskip

\subsection{Translation operator in $L^q_{w}$}

First of all we have to apply Theorem \ref{T1} together with Theorem \ref{TN} to exceptional Laguerre polynomials of the first kind. To this we have to prove that $r(x,t)>0$ on $x>t>0$ for some $m$ and $\alpha$. Recalling that $S(x)=L_m^{(\alpha-1)}(-x)$, we have that $S$ has $m$ simple roots on the negative axis. For simplicity let us denote by $-\xi_i:=x_i$ here, that is $-\xi_m< \dots -\xi_1$. With this notation the following lemma gives a condition for positivity of $r$. After the lemma we give some examples.

\begin{lemma}  In I-Laguerre case $r(x,t)$ (cf. \eqref{34}) is positive on the triangle $\{(x,t)\in\mathbb{R}^2: 0<t<x\}$, if
\begin{equation}\label{F} 1+4(4(m-1)-9+2\alpha)\sum_{i=1}^m\frac{1}{(x+\xi_i)^2}-8\sum_{i=1}^m\frac{\xi_i}{(x+\xi_i)^2}+16\sum_{i=1}^m\frac{\xi_i}{(x+\xi_i)^3}>0.\end{equation}\end{lemma}

\medskip

\proof
Let us recall that $r(x)=x^2+4(\alpha-1+x^2)\frac{S'}{S}(x^2)-4x^2\frac{S''}{S}(x^2)+8x^2\left(\frac{S'}{S}\right)^2(x^2)$. That is it is enough to prove that
$$f(x)=x+4(\alpha-1+x)\frac{S'}{S}(x)-4x\frac{S''}{S}(x)+8x\left(\frac{S'}{S}\right)^2(x)$$
is increasing. Using the differential equation (\cite{ahor})
\begin{equation}\label{des}xS''(x)+(\alpha+x)S'(x)-mS(x)=0,\end{equation}
we have
$$f(x)=x+4(2\alpha-1+2x)\frac{S'}{S}(x)+8x\left(\frac{S'}{S}\right)^2(x)-4m=g(x)-4m.$$
So it it enough to prove that $g(x)$ is increasing.
Since $\frac{S'}{S}(x)=\sum_{i=1}^m\frac{1}{x+\xi_i}$,
$$g'(x)=1+8\sum_{i=1}^m\frac{1}{x+\xi_i}-4(2\alpha-1+2x)\sum_{i=1}^m\frac{1}{(x+\xi_i)^2}+8\sum_{i=1}^m\frac{1}{(x+\xi_i)^2}$$ $$+16\sum_{1\le i<j\le m}\frac{1}{(x+\xi_i)(x+\xi_j)}-16x\sum_{i=1}^m\frac{1}{(x+\xi_i)^3}-16x\sum_{1\le i<j\le m}\frac{2x+\xi_i+\xi_j}{(x+\xi_i)^2(x+\xi_j)^2}.$$
Inserting $\xi_i$ and breaking into partial fractions we have
$$g'(x)=1+8\sum_{i=1}^m\frac{1}{x+\xi_i}-4(2\alpha-1)\sum_{i=1}^m\frac{1}{(x+\xi_i)^2}-8\sum_{i=1}^m\frac{1}{x+\xi_i}+8\sum_{i=1}^m\frac{\xi_i}{(x+\xi_i)^2}$$ $$+8\sum_{i=1}^m\frac{1}{(x+\xi_i)^2}+16\sum_{1\le i<j\le m}\frac{1}{\xi_j-\xi_i}\left(\frac{1}{x+\xi_i}-\frac{1}{x+\xi_j}\right)-16\sum_{i=1}^m\frac{1}{(x+\xi_i)^2}$$ \begin{equation}\label{g}+16\sum_{i=1}^m\frac{\xi_i}{(x+\xi_i)^3}-16x\sum_{1\le i<j\le m}\left(\frac{1}{(x+\xi_i)(x+\xi_j)^2}+\frac{1}{(x+\xi_i)^2(x+\xi_j)}\right).\end{equation}
The last term of the expression above is
$$-16\sum_{1\le i<j\le m}\left(\frac{1}{(x+\xi_i)^2}+\frac{1}{(x+\xi_j)^2}\right)$$ $$+16\sum_{1\le i<j\le m}\left(\frac{\xi_i}{(x+\xi_i)(x+\xi_j)^2}+\frac{\xi_j}{(x+\xi_i)^2(x+\xi_j)}\right)=\Sigma_1+\Sigma_2.$$
Breaking into partial fractions again
$$\Sigma_2=16\sum_{1\le i<j\le m}\left(\frac{\xi_i}{(\xi_j-\xi_i)^2}\left(\frac{1}{x+\xi_i}-\frac{1}{x+\xi_j}\right)-\frac{\xi_i}{\xi_j-\xi_i}\frac{1}{(x+\xi_i)^2}\right.$$ $$\left.+\frac{\xi_j}{(\xi_i-\xi_j)^2}\left(\frac{1}{x+\xi_j}-\frac{1}{x+\xi_i}\right)-\frac{\xi_j}{\xi_i-\xi_j}\frac{1}{(x+\xi_j)^2}\right)$$ $$=-16\sum_{i=1}^m\frac{1}{(x+\xi_i)^2}\sum_{1\le j \le m\atop j\neq i}\frac{\xi_j}{\xi_i-\xi_j}+16\sum_{1\le i<j\le m}-\frac{1}{\xi_j-\xi_i}\left(\frac{1}{x+\xi_i}-\frac{1}{x+\xi_j}\right).$$
Substituting it into \eqref{g} we have
$$g'(x)=1-(8\alpha+36)\sum_{i=1}^m\frac{1}{(x+\xi_i)^2}+8\sum_{i=1}^m\frac{\xi_i}{(x+\xi_i)^2}$$ $$-16\sum_{i=1}^m\frac{1}{(x+\xi_i)^2}\sum_{1\le j \le m\atop j\neq i}\frac{\xi_j}{\xi_i-\xi_j}+16\sum_{i=1}^m\frac{\xi_i}{(x+\xi_i)^3}.$$
Since
$$\sum_{1\le j \le m\atop j\neq i}\frac{\xi_j}{\xi_i-\xi_j}=-(m-1)+\xi_i\frac{\sum_{1\le j \le m\atop j\neq i}\prod_{1\le k \le m\atop k\neq i,j}(\xi_i-\xi_k)}{\prod_{1\le j \le m\atop j\neq i}(\xi_i-\xi_j)}=-(m-1)+\xi_i\frac{\omega''(\xi_i)}{\omega'(\xi_i)},$$
where $\omega(x)=\prod_{i=1}^m(x-\xi_i)$. Since $-\xi_i$-s are the roots of $S$, $\omega(x)=cL_m^{(\alpha-1)}(x)$, the $m$-th  classical Laguerre polynomial. Thus by \eqref{27}
$$\sum_{1\le j \le m\atop j\neq i}\frac{\xi_j}{\xi_i-\xi_j}=-(m-1+\alpha-\xi_i).$$
Finally
$$g'(x)=1+4(4(m-1)-9+2\alpha)\sum_{i=1}^m\frac{1}{(x+\xi_i)^2}-8\sum_{i=1}^m\frac{\xi_i}{(x+\xi_i)^2}+16\sum_{i=1}^m\frac{\xi_i}{(x+\xi_i)^3}.$$

\noindent {\bf Examples.}

\noindent (1) Let $m=1$. Then $S(x)=x+\alpha$, and $g(x)=x+\frac{8\alpha-4+4x}{x+\alpha}+\frac{8x}{(x+\alpha)^2}$.
$$g'(x)=1-\frac{4(1+\alpha)}{(x+\alpha)^2}+\frac{16\alpha}{(x+\alpha)^3}=\frac{1}{(x+\alpha)^3}h(x).$$
Since $h'(x)>3\alpha^2-4\alpha-4$ which is positive when $\alpha>2$, $g$ is increasing for $\alpha>2$ and so $r(x,t)$ is positive in this case.

\noindent (2) Let $m=2$. Then $\frac{S'}{S}(x)=\frac{1}{x+a}+\frac{1}{x+b}$.
$$g'(x)=1+\frac{4}{x+a}+\frac{4}{x+b}-(8\alpha-4+4x)\left(\frac{1}{(x+a)^2}+\frac{1}{(x+b)^2}\right)$$ $$+8\left(\frac{1}{(x+a)^2}+\frac{1}{(x+b)^2}+\frac{2}{(x+a)(x+b)}\right)$$ $$-16x\left(\frac{1}{(x+a)^3}+\frac{1}{(x+b)^3}+\frac{2x+a+b}{(x+a)^2(x+b)^2}\right)$$
Considering that $a=\alpha+1-\sqrt{\alpha+1}$, $b=\alpha+1+\sqrt{\alpha+1}$, we have
$$g'(x)=1-4\frac{a+1}{(x+a)^2}-4\frac{b+1}{(x+b)^2}+16\frac{a}{(x+a)^3}+16\frac{b}{(x+b)^3}.$$
So it is clear that if $\alpha$ is large enough ($\alpha>24$, say) then $r(x,t)$ will be positive again.

\noindent (3) According to \cite[6.31.11]{sz}, if $m>4$ then $\xi_m<4(m-1)+2\alpha$. Since all the terms are decreasing in $\xi_i$
$$\sum_{i=1}^m\frac{4(m-1)+2\alpha-9-2\xi_i}{(x+\xi_i)^2}>m\frac{4(m-1)+2\alpha-9-2\xi_m}{(x+\xi_m)^2}>-m\frac{4(m-1)+2\alpha+9}{(4(m-1)+2\alpha)^2}.$$
$$\frac{4m}{4(m-1)+2\alpha}+\frac{36m}{(4(m-1)+2\alpha)^2}<1,$$
if $\alpha>3m$, say. Then \eqref{F} is positive.

Of course these estimations are rather rough. Our aim was to show that there are several examples for positive $r(x,t)$-s in exceptional Laguerre case.

\medskip

We are in position to state the theorem about exceptional Laguerre translation on $\mathcal{P}$.

\medskip

\begin{corollary}\label{co}Let $u(x,t)=\sum_{k=1}^nu_k(x)u_k(t)$ is the solution on the right upper quadrant to the Cauchy problem
$$u_{xx}(x,t)-u_{tt}(x,t)+\frac{2\alpha+1}{x}u_x(x,t)-\frac{2\alpha+1}{t}u_t(x,t)-r(x,t)=0;$$ $$ \ws\ws u(x,0)=u^0(x)=\sum_{k=1}^nu_k(x); \ws u_t(x,0)=0,$$
where $u(x)=c_{n}L_{m,m+n}^{I, (\alpha)}(x^2)\frac{e^{-\frac{x^2}{2}}}{S(x^2)}$ (cf.\eqref{32}) and $r(x,t)$ is given by \eqref{34}. If
\eqref{F} fulfils, then
$$\|u(x,t)\|_{\infty, [0,\infty]\times [0,\infty]}=\|u^0(x)\|_{\infty, [0,\infty]}.$$\end{corollary}

\medskip

Summarizing the previous ideas we can define the exceptional Laguerre translation via several steps.

\noindent {\it Step 1.} Let $p(x)v(x)=\sum_{k=0}^n a_kR_{m,m+k}^{I, (\alpha)}(x)v(x) \in \mathcal{L}^{I, (\alpha)}_nv$ (for $R_{m,m+k}^{I, (\alpha)}$ see \eqref{n}).\\ Let $u_k(x)=\frac{R_{m,m+k}^{I, (\alpha)}}{S}(x^2)e^{-\frac{x^2}{2}}$.\\ Then $u(x,t)=\sum_{k=0}^n a_ku_k(x)u_k(t)$ is the symmetric solution of
$$u_{xx}-u_{tt}+\frac{2\alpha+1}{x}u_x - \frac{2\alpha+1}{t}u_t-ru=0,$$
$$ u(x,0)=(pv)(x^2)=\sum_{k=0}^n a_ku_k(x), \ws u_t(x,0)=0$$
on $[0,\infty)$, where $r(x,t)$ is given in \eqref{34}. That is $T_t(pv,x^2)=(Tpv)(x^2,t^2)=u(x,t)$ which defines the translation operator on the real linear space $\mathcal{P}$ generated by $\{u_k(x)\}_{k}$ on $[0,\infty)$. By Corollary \ref{co} the operator norm is one on $\mathcal{P}$.

\noindent {\it Step 2.} Since $x,t\ge 0$ by the equivalence $x^2 \to x$, $t^2 \to t$ this operator uniquely defines an operator - denoting by $T_t$ again - which acts on $\mathcal{L}^{I, (\alpha)}v$ and maps $(pv)(x)$ to $(pv)_t(x)=(pv)(x,t)=\sum_{k=0}^n a_kR_{m,m+k}^{I, (\alpha)}(x)v(x)R_{m,m+k}^{I, (\alpha)}(t)v(t)\in C_1[0,\infty)$. The operator norm is at most one for all $t\ge 0$ again.

\noindent {\it Step 3.} According to Lemma \ref{L3} $\mathcal{L}^{I, (\alpha)}v$ is dense in $C_1[0,\infty)$, that is the space of continuous functions disappear at infinity. Thus the extension of the operator from $\mathcal{L}^{I, (\alpha)}v$ to $C_1[0,\infty)$ can be done with the same norm.

\noindent {\it Step 4.} Let us denote by
$$\tilde{f}(x)=(fv)(x).$$
Let $p,q \in \mathcal{L}^{I, (\alpha)}$ arbitrary and $n:=\max\{\deg p,\deg q\}$. By orthogonality
$$\langle T_t(\tilde{p},x)\tilde{q}(x)\rangle_w=\int_0^{\infty}\sum_{k=0}^na_kR_{m,m+k}^{I, (\alpha)}(x)\tilde{R}_{m,m+k}^{I, (\alpha)}(t)\sum_{j=0}^nb_jR_{m,m+j}^{I, (\alpha)}(x)x^{\alpha}\frac{e^{-x}}{S^2(x)}dx$$ $$=\int_0^{\infty}\sum_{k=0}^na_kR_{m,m+k}^{I, (\alpha)}(x)\sum_{j=0}^nb_jR_{m,m+j}^{I, (\alpha)}(x)\tilde{R}_{m,m+j}^{I, (\alpha)}(t)x^{\alpha}\frac{e^{-x}}{S^2(x)}dx=\langle \tilde{p}(x)T_t(\tilde{q},x)\rangle_w,$$
cf. \eqref{s}. By duality, similarly to \eqref{d} we show that for all $t\ge 0$ $\|T_t\|_{1,w}\le 1$. Since for all $\varepsilon>0$ there is a $g\in L_{[0,\infty),\infty}$ such that $\langle T_t(\tilde{p},x)g(x)\rangle_w>\|T_t(\tilde{p},x)\|_{w,1}-\frac{\varepsilon}{2}$, and for this $\varepsilon>0$ there is an  $M>0$ such that $\langle T_t(\tilde{p},x)g(x)\rangle_w-\frac{\varepsilon}{2}<\int_0^M T_t(\tilde{p},x)g(x)w(x)dx$ and again for this $\varepsilon>0$ there is a $q\in \mathcal{L}^{I, (\alpha)}$ such that $\left|\int_0^M T_t(\tilde{p},x)g(x)w(x)dx-\int_0^{\infty} T_t(\tilde{p},x)\tilde{q}(x)w(x)dx \right|<\varepsilon$, we can write that
$$\|T_t(\tilde{p},x)\|_{1,w}=\sup_{q\in \mathcal{L}^{I, (\alpha)} \atop \|\tilde{q}\|_{\infty}=1}\langle T_t(\tilde{p},x)\tilde{q}(x)\rangle_w$$ $$=\sup_{q\in \mathcal{L}^{I, (\alpha)} \atop \|\tilde{q}\|_{\infty}=1}\langle \tilde{p}(x)T_t(\tilde{q},x)\rangle_w\leq \sup_{q\in \mathcal{L}^{I, (\alpha)} \atop \|\tilde{q}\|_{\infty}=1}\|T_t\tilde{q}\|_{\infty}\| \tilde{p}\|_{1,w}\leq \| \tilde{p}\|_{1,w}.$$
After these by density we can extend the operator with the same norm to $L^1_w$.

\noindent {\it Step 5.} Now to keep the $p$-norm of the operator under one we want to apply the complex version of the Riesz-Thorin theorem. So we have to extend the operator to complex-valued functions such that the extended operator is complex linear, and the infinity- and one-norm of the operator does not grow. So let $f: [0,\infty) \to \mathbb{C}$, $f(x)=u(x)+iv(x)$ such that for all $t\ge 0$ $T_t(u,x)=u(x,t)$, then the extended operator - denoted by $T_t$ again is $T_t(f,x)=u(x,t)$. It is a real-valued additive and real homogeneous operator and since $|u(x)|\leq |f(x)|$, the operator norms are unchanged.
Below we proceed as in \cite[Lemma 1]{ar}. Let $\hat{T_t}f:=T_tf-iT_t(if)$. Then $\hat{T_t}$ is an extension of $T_t$ and it can be easily seen that it is linear with respect to complex coefficients. For a fixed $x\in [0,\infty)$ $\hat{T_t}(f,x)=re^{i\varphi}$, that is since $\hat{T_t}$ is complex homogeneous and $T_t$ is real-valued,
$$|\hat{T_t}(f,x)|=e^{-i\varphi}\hat{T_t}(f,x)=T_t(e^{-i\varphi}f,x)-iT_t(ie^{-i\varphi}f,x)=T_t(e^{-i\varphi}f,x).$$

The infinite-norm of the complex operator restricted to $C_1$- functions is also one. Indeed, let $f$ be such that for $\Re f=u, \Im f=v\in C_1$, that is $\lim_{x\to\infty}u(x)=\lim_{x\to\infty}v(x)=0$. Since by density
$\|T_t|_{C_1}\|_{\infty}=1$, there is a $\xi\in [0,\infty)$ such that
$$|\hat{T_t}(f,x)|=T_t(e^{-i\varphi}f,x)=T_t(\cos\varphi u+\sin \varphi v,x)\le |(\cos\varphi u)(\xi)+(\sin\varphi v)(\xi)|$$ $$\le \left((\cos\varphi u)(\xi)+(\sin\varphi v)(\xi))^2+(\cos\varphi v)(\xi)-(\sin\varphi u)(\xi))^2\right)^{\frac{1}{2}}\leq \|f\|_{\infty},$$
that is $\|\hat{T}_t|_{C_1}\|_{\infty}=1$.

Investigating the $1$-norm let $g(x):=g_t(x)=e^{-i\varphi(x)}f(x)$, where $\varphi(x)$ is the argument of $\hat{T_t}(f,x)$ with a fixed $t$. Then considering that $\|T_t\|_{1,w}=1$
$$\int_0^{\infty}|\hat{T_t}(f,x)|w(x)dx=\int_0^{\infty}|T_t(g,x)|w(x)dx\le\int_0^{\infty}|g(x)|w(x)dx=\|f\|_{1,w},$$
that is $\|\hat{T}_t\|_{1,w}=1$.

\noindent By the same procedure we get that  $\|\hat{T}_t\|_{2,w}=1$.

\noindent {\it Step 6.} Now we are in position to apply the complex Riesz-Thorin theorem between $L^1_w$ and $L^2_w$, and then we get that  $\|\hat{T_t}\|_{q,w}=1$ for all $1\le q\le 2$ and $t\ge 0$, thus the same is valid for $T_t$.

\noindent {\it Step 7.} The Riesz-Thorin theorem cannot be applied automatically between $L^2$ and $L^{\infty}$, since the operator norm is given only on a subspace of $L^{\infty}$. So we repeat the computation above, that is let $2<q<\infty$ and $\frac{1}{q}+\frac{1}{q'}=1$. By Lemma \ref{L3}
$$\|T_t(\tilde{p},x)\|_{q,w}=\sup_{q\in \mathcal{L}^{I, (\alpha)} \atop \|\tilde{q}\|_{q'}=1}\langle T_t(\tilde{p},x)\tilde{q}(x)\rangle_w$$ $$=\sup_{q\in \mathcal{L}^{I, (\alpha)} \atop \|\tilde{q}\|_{q'}=1}\langle \tilde{p}(x)T_t(\tilde{q},x)\rangle_w\leq \sup_{q\in \mathcal{L}^{I, (\alpha)} \atop \|\tilde{q}\|_{q'}=1}\|T_t\tilde{q}\|_{q',w}\| \tilde{p}\|_{w,q}\leq \| \tilde{p}\|_{q,w}.$$
Now we can extend the operator from polynomials to the whole $L^q_w$ space with the same norm, as above. Summarizing we have

\begin{theorem} The exceptional Laguerre translation operator can be extended from exceptional polynomials to the weighted $L^q_w$ spaces with norm one, that is
$$\|T_t\|_{q,w}=1, \ws \ws1\le q <\infty$$ and
$$\|T_t|_{C_1}\|_{\infty}=1.$$\end{theorem}

\subsection{Nikol'skii inequality}
In this subsection similarly to  the investigations of \cite{ad}, \cite{ad2} and \cite{adh} we characterize the extremal elements in the endpoint- and norm Nikol'skii-type inequalities.\\ The problems are as it follows. Let - as a linear space over $\mathbb{C}$
$$\mathcal{L}_{n,\mathbb{C}}^{I, (\alpha)}:=\mathrm{span}\left\{L_{m,m+k}^{I, (\alpha)}(x)\right\}_{k=0}^n, \ws \ws \mathcal{L}_{\mathbb{C}}^{I, (\alpha)}=\cup_n\mathcal{L}_{n,\mathbb{C}}^{I, (\alpha)}.$$

\noindent {\bf (P1)}: We look for the best constant in the inequality
\begin{equation}\label{p1} |\tilde{p}(x)|\le D_{n,q}(x)\|\tilde{p}\|_{q,w}, \end{equation}
where $p\in \mathcal{P}_n$.

\noindent {\bf(P2)}: We look for the best constant in the inequality
\begin{equation}\label{p2} \|\tilde{p}\|_{\infty}\le M_{n,q}\|\tilde{p}\|_{q,w}, \end{equation}
where $p\in  \mathcal{P}_n$.

Let us denote by
$$\mathcal{P}_n^{(z)}:=\left\{p\in \mathcal{P}_n : p(z)=0\right\}.$$

The connection between (P1) and (P2) is given by the following application of the beautiful theorem of Arestov, \cite[Theorem 1]{ar}:

\begin{theorem}\label{a}\cite[Theorem 4]{ar} Let $1\le q <\infty$.  $\tilde{\varrho}_n$ is extremal in $(\mathrm{P}1)$ if and only if
$$\int_0^{\infty}|\tilde{\varrho}_n(x)|^{q-1}\mathrm{sign}\varrho_n(x)p(x)w(x)dx=0, \ws \ws \forall \ws p\in \mathcal{P}_n^{(z)}.$$
If $1<q <\infty$, the extremal element is unique (up to a constant factor).\end{theorem}

\remark

\noindent (1) $\frac{1}{D_{n,q}(x)}$ is just the so-called $q$-Christoffel function (cf. e.g. \cite{n}), i.e. for $q=2$ it is the Cristoffel function, which is the reciprocal of the square-sum of the first $n$ orthogonal polynomials in classical and generalized cases. The Christoffel function is also in close connection (the $x=t$ case) with the kernel function of Fourier sums which can be investigated by Christoffel-Darboux formula for classical and generalized orthogonal polynomials. The investigation of $D_{n,q}(x)$ might be important from this point of view, since Christoffel-Darboux formula is a consequence of the three-term recurrence formula, which does not exist in case of exceptional orthogonal polynomials. (In our case there is a five-term recurrence formula.)

\noindent (2) Here we restrict our investigation to $D_{n,q}(0)$. Zero, as the endpoint of our interval plays a special role. The assumption $u_k(0)=1$ for all $k$ ensures that the initial condition fulfils (cf. \eqref{i}). The assumption $\|u_k\|_{\infty}=1$ ensures that the operator norm is one as it is derived in $L^2$ case in finite dimension (cf. \eqref{n1}). It means that each of the elements of the orthogonal system attains its maximum at the same point, cf. Lemma \ref{Ll1} and Lemma \ref{Ll2}. This special point is zero. We mention here that at least in cases of classical and generalized orthogonal polynomials it can happen only at the endpoint:

Indeed, let $u_n(x)=p_nw(x)$, where $w>0$ on $(x_0,x_1)$ and $\{p_n\}$ fulfils the three-term recurrence formula. Let us assume that this common maximum/minimum  point $t_0$ is in the interior of the interval. Then
$$\left(wp_n\right)'(t_0)=\left(w((Ax+B)p_{n-1}+Cp_{n-2})\right)'(t_0)$$ $$=A(wp_{n-1})(t_0)+(Ax+B)\left(wp_{n-1}\right)'(t_0)+C\left(wp_{n-2}\right)'(t_0)=0,$$
where $A,B,C$ are the coefficients in the recurrence formula. Since  $t_0$ was a common extremum point and $w(t_0)>0$,
$$p_{n-1}(t_0)=0.$$
On the other hand
$$0=(wp_{n-1})'(t_0)=(w'p_{n-1})(t_0)+(wp'_{n-1})(t_0),$$
that is
$$p'_{n-1}(t_0)=0,$$ which is impossible since the orthogonal polynomials have only simple zeros.

\medskip

After this we turn to an application of exceptional Laguerre translation.

\begin{theorem}\label{m} For all $1\le q <\infty$,
$$D_{n,q}(0)=M_{n,q},$$
the extremal polynomial both in the inequalities \eqref{p1} and \eqref{p2} is $\tilde{\varrho}_{n}$ which for $q>1$ is unique and attains its maximum at zero.
\end{theorem}

\proof
Let $p_n \in \mathcal{L}^{I, (\alpha)}_{n,\mathbb{C}}$ arbitrary. Using that $\|T_t|_{\mathcal{P}}\|_{\infty}=1$ and applying standard arguments, by Lemma \ref{Ll1} and Corollary \ref{co}
$$\|\tilde{p_n}\|_{\infty}=|\tilde{p_n}(t)|=\left|\sum_{k=0}^na_k\tilde{R}_{m,m+k}^{I,\alpha}(t)\tilde{R}_{m,m+k}^{I,\alpha}(0)\right|=|T_t(\tilde{p_n},0)|=|\tilde{r_n}(0)|$$ $$\le D_{n,q}\|\tilde{r}_n\|_{q,w}\le  D_{n,q}\|T_{t}\|_q\|\tilde{p}_n\|_{q,w}= D_{n,q}\|\tilde{p}_n\|_{q,w}.$$
That is $ M_{n,q}\le D_{n,q}$. (The other direction is obvious.) Applying this to $\tilde{\varrho}_{n}$ (cf. Theorem \ref{a}) all the inequalities become equalities, thus for $q>1$ $p=\varrho_n=r$, $t=0$ and $\tilde{\varrho}_{n}$ attains its maximum at zero.

\medskip

\remark
Theorem \ref{a} shows that the extremal polynomial coincides with the Chebyshev polynomial in the weighted space in question, i.e. the polynomial with leading coefficient $1$ which deviates from zero the least, cf. e.g. \cite{ar}, \cite{adh}.

\medskip

\noindent \small{Department of Analysis, \newline
Budapest University of Technology and Economics}\newline
\small{ g.horvath.agota@renyi.mta.hu}

\end{document}